# GROUPS IN SIMPLE THEORIES

Frank O. Wagner

ABSTRACT. Groups definable in simple theories retain the chain conditions and decomposition properties known from stable groups, up to commensurability. In the small case, if a generic type of $G$ is not foreign to some type $q$, there is a $q$-internal quotient. In the supersimple case, the Berline-Lascar decomposition works. One-based simple groups are finite-by-abelian-by-finite.

## INTRODUCTION

In connection with the work by Hrushovski and Chatzidakis [CH1] on difference fields, and Cherlin and Hrushovski [CH2] on smoothly approximated structures, simple theories of groups have recently begun to attract attention. While the above concentrate on those particular important examples, all of which have finite rank, simple theories in general have been studied by Shelah [Sh1,2], Kim and Pillay [K1-3,KP], and in the group case by Pillay [Pi].

In this paper, we shall reconstitute the basic theory (chain conditions, internal quotients, components, Zil'ber Indecomposability, one-based groups) for groups with simple theories, some under the additional condition of supersimplicity or smallness. We shall freely use notions and results from the papers by Kim and Pillay. In particular, we shall make use of the local ranks $D(.,\varphi,k)$ and the theory of generic types and stabilizers, as developed in [Pi].

## 1. CHAIN CONDITIONS

A family $\mathfrak{F}$ of subgroups of a group $G$ is said to be *uniformly definable* if every element $F_i$ of $\mathfrak{F}$ is of the form $\{g \in G : G \models \varphi(g, \bar{a}_i)\}$, where the formula $\varphi$ remains fixed and only the parameters $\bar{a}_i$ vary. The most general chain condition satisfied by a stable group is the chain condition on uniformly definable subgroups:

[**ucc**] Any chain of uniformly definable subgroups has finite length, bounded by some $n$ depending only on the defining formula $\varphi$.

Clearly, the formula $\forall x\ [\varphi(x,\bar{y}) \to \varphi(x,\bar{z})]$ defines a partial order; as a simple theory cannot interpret a partial order with an infinite chain, the ucc must hold in any simple theory as well. However, the more useful chain condition in a stable theory is the chain condition on intersections of uniformly definable subgroups:

[**icc**] Any chain of intersections of uniformly definable subgroups has finite length, bounded by some $n$ depending only on the defining formula $\varphi$.

This need not hold in a simple group, as the following example shows:

**Example.** Let $A = C_p^I$ be an elementary abelian group of exponent $p$, and consider a predicate $P(x, j) = \{(a_i)_{i \in I} \in A : a_j = 0\} \subset A \times I$. The structure $\langle A \cup I, 0, +, P \rangle$ is simple, but has an infinite descending chain of intersections of uniformy definable subgroups.

However, $A$ obviously has no infinite descending chain of definable subgroups, each of infinite index in the previous one. In general, we shall have to consider type-definable subgroups as well, in order to recover the icc up to commensurability.

**Definition 1.1.**
  (1) A group is *type-definable* if it is given as the set of realizations of a partial type (in any elementary superstucture). It is $\bigcap$-definable if it is given as the intersection of definable groups.
  (2) If $H \leq G$ are two type-definable groups, then $H$ has *bounded index* in $G$ if for any elementary superstructure $\mathcal{M}$ of the ambient theory the index $|G^{\mathcal{M}} : H^{\mathcal{M}}|$ remains bounded.

Clearly a $\bigcap$-definable group is type-definable. In a stable theory, the converse holds as well; this is open in general for simple theories. Note that the infinitesimals form a type-definable subgroup of $\langle \mathbb{Q}, 0, +, < \rangle$ which is not $\bigcap$-definable (but this structure is not simple). If $H$ is $\bigcap$-definable as $\bigcap_{i < \kappa} H_i$ and has bounded index in $G$, then by compactness every $H_i$ intersects $G$ in a subgroup of finite index. So the index of $H$ in $G$ is bounded by $2^\kappa$. This bound still holds if $H$ is type-defined by a partial type $\pi$ of size $\kappa$: if $\varphi(x) \in \pi$, then there is some definable superset $X$ of $G$ and some $n < \omega$ such that for any $n$ elements $g_1, \ldots, g_n \in X$ we have $\models \varphi(g_i g_j^{-1})$ for at least one pair $i \neq j$, as otherwise we could find models in which $H$ has arbitrarily large index in $G$ by compactness. Now suppose $\{g_i : i \in (2^\kappa)^+\}$ were a set of representatives for different cosets of $H$ in $G$. Then for any pair $\{i, j\}$ we could associate a formula $\varphi \in \pi$ such that $\models \neg\varphi(g_i g_j^{-1})$. By the Erdös-Rado Theorem [ER] there is a subset $J \subset (2^\kappa)^+$ of cardinality $\kappa^+$ such that all pairs from $J$ are associated with the same formula. This contradicts the fact that we even have a finite bound on such sets.

**Fact 1.2.** *[Pi, 3.12] Suppose $H \leq G$ are type-definable subgroups in a simple theory. Then the index of $H$ in $G$ is bounded iff $D(H, \varphi, k) = D(G, \varphi, k)$ for all formulas $\varphi$ and all parameters $k < \omega$.*

**Proposition 1.3.** *Let $(H_i : i \leq \alpha)$ be a descending chain of type-definable subgroups in a simple group, continuous at limits, such that each successor group has unbounded index its predecessor. Then $\alpha < |T|^+$.*

*Proof.* For any $i < \alpha$ let $(\varphi_i, k_i)$ be a pair such that $D(H_i, \varphi_i, k_i) > D(H_{i+1}, \varphi_i, k_i)$. If the chain had length $|T|^+$, then we could find a subchain $J$ of length $|T|^+$ such that $(\varphi_i, k_i) = (\varphi, k)$ is constant for $i \in J$. As $D(H_j, \varphi, k)$ is strictly descending for $j \in J$, this contradicts finiteness of $D(T, \varphi, k)$. □

For $\bigcap$-definable groups, we get stronger results.

**Definition 1.4.** Two groups $G$ and $H$ are commensurable if $G \cap H$ has finite index in both $G$ and $H$. A group $G$ is uniformly commensurable to a family $\mathfrak{F}$ if the index of $G \cap H$ in $G$ and in $H$ is bounded independently from $H \in \mathfrak{F}$.

We shall need a theorem due to Schlichting and Bergman-Lenstra:

**Fact 1.5.** *Let $G$ be a group and $\mathfrak{F}$ a family of subgroups such that there is $n < \omega$ bounding the index $|F : F \cap F'|$ for any $F, F' \in \mathfrak{F}$. Then there is a subgroup $N$ which is uniformly commensurable to $\mathfrak{F}$ and invariant under all automorphisms of $G$ which stabilize $\mathfrak{F}$ setwise.*

In fact, it is easily seen from the proof (or from commensurability) that $N$ is definable, and contained in the product $FF'$ for two elements $F, F' \in \mathfrak{F}$.

**Lemma 1.6.** *Let $G_0$ be a definable group in a simple theory, $G$ a type-definable subgroup, and $H(a)$ an $a$-definable subgroup. Then there is some integer $n < \omega$ and a definable superset $X$ of $G$ such that whenever $H(a')$ intersects $G$ in a subgroup of finite index, the index is bounded by $n$ and $X$ is covered by $n$ cosets of $H(a')$.*

*Proof.* The condition on $a'$ that $H(a')$ intersects $G$ in a subgroup of finite index is both open (finitely many translates of $H(a')$ cover some definable superset of $G$) and closed ($D(G \cap H(a'), \varphi, k) = D(G, \varphi, k)$ for all $\varphi, k$). By compactness, it is definable, and there is a bound $n$ and a definable $X \supseteq G$ as required. □

**Proposition 1.7.** *Let $G$ be a definable group in a simple theory, and $\mathfrak{F}$ a family of uniformly definable subgroups. Then there is an definable subgroup $N$ and a finite intersection $N_0 := \bigcap_{i<n} H_i$ of elements in $\mathfrak{F}$, such that $N$ is uniformly commensurable with $H \cap N_0$ for all $H \in \mathfrak{F}$, and $N$ is invariant under all automorphisms stabilizing ($G$ and) $\mathfrak{F}$.*

In particular, any intersection of a family of uniformly definable subgroups is finite up to commensurability (where commensurability for type-definable groups means that the intersection has *bounded* rather than *finite* index in either group).

*Proof.* By Proposition 1.3 there is an intersection $K$ of groups in $\mathfrak{F}$ of size at most $|T|$, such that any bigger intersection of this form yields a subgroup of bounded index. Hence for any $H \in \mathfrak{F}$ the index $|K : K \cap H|$ is finite, and by Lemma 1.6 there is $k < \omega$ and a finite subintersection $N_0 = \bigcap_{i<n} H_i$ of $K$ such that $|N_0 : N_0 \cap H| \leq k$ for all $H \in \mathfrak{F}$.

Now consider the family of conjugates of $N_0$ under automorphisms stabilizing $\mathfrak{F}$. For any such conjugate $N_0'$, the index $|N_0 : N_0 \cap N_0'|$ is bounded by $k^n$. By Fact 1.5, there is a definable subgroup $N$ uniformly commensurable with every element in the family, and invariant under all these automorphisms. □

In particular, if $G$ is $|T|^+$-saturated and $\mathfrak{F}$ is invariant under all automorphisms fixing some parameter set $A$, then $N$ is $A$-definable.

**Corollary 1.8.** *Let $\mathfrak{F}$ be the family of automorphic conjugates of some definable subgroup $H$ in a simple theory. Then the equivalence relation on $\mathfrak{F}$ given by commensurability is definable. Furthermore, there is a subgroup $H^c$ commensurable with $H$ such that any conjugate of $H^c$ is either equal to $H^c$ or intersects it in a subgroup of infinite index.*

Hence $H^c$ may be considered some kind of locally connected component of $H$ — however, it need not be unique, and is not necessarily a finite intersection of commensurable conjugates of $H$.

*Proof.* Let $\mathfrak{F}_0$ be the subfamily of automorphic conjugates of $H$ commensurable with $H$. Then the index of $H \cap H'$ in $H$ is bounded by some $n < \omega$ for all $H' \in \mathfrak{F}_0$ by Lemma 1.6, so commensurability is definable. Furthermore, Fact 1.5 yields a definable subgroup $H^c$ commensurable with $H$, such that $H^c$ is stabilized by any automorphism stabilizing $\mathfrak{F}_0$. In particular, any commensurable conjugate of $H^c$ is equal to $H^c$, as it must arise from the same family $\mathfrak{F}_0$. $\square$

As a consequence, we recover the icc up to commensurability.

**Theorem 1.9.** *Let $G$ be a definable group in a simple theory, and $\mathfrak{F}$ a family of uniformly definable subgroups of $G$. Then there are integers $k, k' < \omega$ such that any intersection $N$ of elements in $\mathfrak{F}$ is commensurable with a subintersection $N_0$ of size at most $k$ and $|N_0 : N_0 \cap H| \leq k'$ for any $H \in \mathfrak{F}$ with $H \geq N$. Furthermore, there is an integer $n < \omega$ such that any chain of intersections of elements of $\mathfrak{F}$, each of infinite index in its predecessor, has length at most $n$.*

*Proof.* Clearly, we may assume that $\mathfrak{F}$ is a maximal family, i.e. the family of all subgroups definable by a particular formula. By Proposition 1.7 every intersection $N$ of elements of $\mathfrak{F}$ is commensurable with a finite subintersection $N_0$. Suppose there is no bound $k$ on the size of this subintersection, or no bound $k'$ on the index $|N_0 : N_0 \cap H|$ for $H \in \mathfrak{F}$ with $H \geq N$. Then by compactness, the following type is consistent (where we identify a group $f \in \mathfrak{F}$ with its parameter):

$$\{F_i \in \mathfrak{F} : i \in \omega\} \cup \{|\bigcap_{i<k} F_i : \bigcap_{i \leq k} F_i| > k' : k, k' \in \omega\}.$$

But any realization of this type yields a family which contradicts Proposition 1.7. So there are the required bounds.

Now consider a chain $C_0 > C_1 > \cdots > C_n$ of intersections of groups in $\mathfrak{F}$ such that every $C_i$ has infinite index in its predecessor. By the first part, we may replace every $C_i$ with a uniformly definable and uniformly commensurable $C'_i$. These $C'_i$ may no longer form a chain, but the index $|C'_i : C'_i \cap C'_j|$ for $j \leq i$ is bounded by some $k_0 = k^{k'}$ which only depends on the initial family $\mathfrak{F}$. Suppose these $C'_i$ are given by some formula $\varphi(x, a_i)$. Consider the relation

$$a < a' \text{ iff } |\varphi(x,a) : \varphi(x,a) \cap \varphi(x,a')| \leq k_0 \wedge |\varphi(x,a') : \varphi(x,a) \cap \varphi(x,a')| > k_0.$$

Clearly $<$ does not have any triangles $a < b < c < a$; since $C_j$ has infinite index in $C_i$ for $i < j$, we have $a_i > a_j$ iff $i < j$. If $n$ were not bounded, this would contradict simplicity by [Sh2, Claim 2.7]. It follows that there must be a bound on the length of the initial chain $(C_i : i \leq n)$. $\square$

Thus, in a simple theory we may recover most of the components known from stable groups, but only up to commensurability.

We finish with two general lemmas, which are more or less implicit in [CH2]:

**Lemma 1.10.** *Suppose $G$ is a type-definable group in a simple theory, and $X$ is a type-definable subset of $G$ such that for independent $x, x' \in X$ we have $x^{-1}x' \in X$.*

*Then $X \cdot X =: Y$ is a type-definable subgroup of $G$, and $X$ is generic in $Y$. In fact, $X$ contains all generic types for $Y$.*

*Proof.* As $X' = \{x \in X : x^{-1} \in X\}$ also satisfies the assumptions of the lemma and generates the same group (if $a, b, c$ are three independent elements of $X$, then $a^{-1}b \in X'$ and $(b^{-1}a)c \in X'$, whence $a \in X' \cdot X'$), we may assume that $X$ is closed under inversion.

Enumerate all pairs of stratified formulas and natural numbers as $(\varphi_i, k_i)$, for $i < \kappa$. Choose a type $p$ containing the formulas $x \in X$ such that $D(p, \varphi_i, k_i) = n_i$ is maximal, subject to $d(p, \varphi_j, k_j) = n_j$ for all $j < i$. Then if $a, b, c$ are three elements of $X$, choose $d \models p$ with $d \downarrow a, b, c$. Then $bd \in X$. As the stratified ranks of $\mathrm{tp}(d/b)$ are those of $\mathrm{tp}(bd/b)$ and are equal to those of $\mathrm{tp}(bd)$ by maximality, we get $bd \downarrow b$, whence $bd \downarrow a, b$ and $abd \in X$. Finally $d \downarrow c$, so $d^{-1}c \in X$, and $abc \in Y$.

Now let $a$ be a generic element of $Y$, say $a = bc$ for some $b, c \in X$. Let $d \models p$ independently from $b, c$. Then $ad$ is again a generic element. On the other hand $cd \in X$, and as above we get that $\mathrm{tp}(cd/c)$ and $\mathrm{tp}(cd)$ have the same stratified ranks, whence $cd \downarrow c$. Therefore $b \downarrow cd$, and $bcd = ad \in X$. But since $a$ is generic independent from $d$, we have $ad \downarrow d$, so $add^{-1} = a \in X$. □

This will yield most of the type-definable groups we shall encounter. In order to normalize them (as we cannot use Fact 1.5 and the subsequent results) we shall rely on the following.

**Lemma 1.11.** *Let $G$ be a group definable in a simple theory, and $H$ a subgroup type-definable over some parameters $A$. Then there is a type-definable group $N$ invariant under strong $\emptyset$-automorphisms (i.e. $N$ is definable over any model) commensurable with sufficiently big intersections of strong $\emptyset$-conjugates of $H$.*

*Proof.* By Proposition 1.3 there is an intersection $N(A)$ of strong $\emptyset$-conjugates of $H$, type-definable over some parameters $A$, such that any further intersection has bounded index in $N(A)$. Let $X$ be the set

$$\{g \in G : g \in N(A') \text{ for some } A' \text{ of the same Lascar strong type as } A \text{ with } g \downarrow A'\}.$$

If $g, g' \in X$ and $g \downarrow g'$, then by the Independence Theorem we may choose $A'$ of the same Lascar strong type as $A$ and independent from $g, g'$ such that both $g$ and $g'$ are in $N(A')$, whence $g^{-1}g' \in N(A')$; since $g^{-1}g' \downarrow A'$, we get $g^{-1}g' \in X$. By Lemma 1.10 the type-definable set $X^2$ forms a subgroup $N$, which is clearly invariant under strong automorphisms. We claim that $N$ is commensurable with $N(A)$.

Consider $A'$ of the same Lascar strong type as $A$ and independent from $A$. By minimality, $N(A') \cap N(A)$ has bounded index in $N(A)$ and therefore is generic for $N(A)$. So there is an element $g \in N(A)$ generic over $AA'$ which lies in $N(A')$. But then $A' \downarrow_A g$ and $A' \downarrow A$ yields $A' \downarrow g$, whence $g \in X$. So $N \cap N(A)$ is generic in $N(A)$. Conversely, let $g$ be generic in $N$ over $A$. Then there is $A'$ of the same Lascar strong type as $A$ with $g \downarrow A'$ and $g \in N(A')$. But then for every $\varphi$ and $k$

$$D(N, \varphi, k) = D(g/\emptyset, \varphi, k) = D(g/A', \varphi, k) \leq D(N(A'), \varphi, k) = D(N(A), \varphi, k).$$

Therefore $D(N, \varphi, k) \leq D(N(A), \varphi, k) = D(N(A) \cap N, \varphi, k)$; equality of the ranks and commensurability follow. □

## 2. Small Groups

Recall that a theory is *small* if it has only countably many pure types. Small stable groups are well-behaved, and this is similar in the simple case.

**Fact 2.1.** *[Ki3] In a small theory every type-definable equivalence relation over a finite set is the conjunction of definable ones.*

**Fact 2.2.** *[HKP] In a simple theory, for every type $p$ over a finite set there is a finite tuple $\bar{a}$ and a type-definable equivalence relation $E$ such that whenever a non-forking extension $q$ of $p$ does not fork over a set $B$, then $E$ has boundedly many classes on $tp(\bar{a}/B)$, and $p$ does not fork over $\bar{a}/E$.*

An equivalence class modulo a type-definable equivalence relation is called a *hyperimaginary*. Fact 2.1 tells us that small theories eliminate hyperimaginaries: every hyperimaginary is really just a set of imaginary elements. Fact 2.2 says that types over finite sets have hyperimaginary canonical bases (in fact arbitrary types have countable hyperimaginary canonical bases); in conjunction with Fact 2.1 this yields the existence of canonical bases for types over finite sets in a small simple theory. Furthermore, Lascar strong types over finite sets are just strong types in a small simple theory.

**Lemma 2.3.** *Suppose the small group $G$ has a subgroup $H$ which is type-definable over a finite set. Then $H$ is $\bigcap$-definable.*

*Proof.* On $G \times G$ the relation $xy^{-1} \in H$ is type-definable, and clearly an equivalence relation. By Fact 2.1, it is the intersection of definable equivalence relations $E_i$. Put
$$H_i = \{g \in G : \forall x\, xE_i gx)\,\}.$$

If $h \in H$, then clearly $(hx)x^{-1} \in H$, whence $xE_i hx$ for all $i$, and $H \subseteq H_i$. On the other hand, if $h \in H_i$ for all $i$, then $1E_i h$ for all $E_i$ (put $x=1$), whence $h \in H$. Furthermore, if $g, h \in H_i$, then $h^{-1}xE_i hh^{-1}x$, and $xE_i hxE_i ghx$. It follows that $H_i$ is a subgroup of $G$, and $H = \bigcap_i H_i$. □

**Definition 2.4.** *Let $Q$ be an $A$-invariant family of types. A type $p \in S(A)$ is $Q$-internal if for every realization $a$ of $p$ there is $B \downarrow_A a$, types $\bar{q}$ from $Q$ based on $A$, and realizations $\bar{c}$ of $\bar{q}$, such that $a \in \mathrm{dcl}(B\bar{c})$.*

Note that $B$ may well depend on the particular realization of $p$, contrary to the usual definition of internality in a stable theory, where $B$ is fixed before choosing $a$ (and the two notions are equivalent). This arises from the fact that in a simple theory a type (even over a model) may have many different non-forking extensions.

**Proposition 2.5.** *Suppose $tp(a)$ is non-orthogonal to a type $q$ in a small simple theory, and let $Q$ be the family of $\emptyset$-conjugates of $q$. Then there is $a_0 \in acl(a) - acl(\emptyset)$ such that $tp(a_0)$ is $Q$-internal.*

*Proof.* This is just the usual. Let $tp(c/A)$ be a non-forking extension of $q$ such that $a \not\downarrow_A c$, and consider a Morley sequence $(c_i A_i : i < \omega)$ in $tp(cA/a)$. Then the canonical base $Cb(cA/a)$ is algebraic in $a$ and definable over $(c_i A_i : i < \omega)$. Since $a \not\downarrow cA$, this canonical base is not contained in $\emptyset$ and there is some finite bit $a_0 \in \mathrm{acl}(a) - \mathrm{acl}(\emptyset)$ definable over a finite tuple $\bar{c}\bar{A}$ of the Morley sequence. Now $a_0 \underset{}{\downarrow} \bar{A}$ and every co-ordinate of $\bar{c}$ realizes an $\emptyset$-conjugate of $q$ over $\bar{A}$. It follows that $\mathrm{tp}(a_0)$ is $Q$-internal. $\square$

**Theorem 2.6.** *Let $G$ be a $\bigcap$-definable group over a finite set $A$ in a small simple theory, and suppose that a generic type $p$ of $G$ is non-orthogonal to some type $q$. Then there is a normal relatively definable subgroup of infinite index in $G$ such that the quotient $G/N$ is $Q$-internal, where $Q$ is the collection of $A$-conjugates of $q$.*

*Proof.* If $Q$ denotes the set of $A$-conjugates of $q$ and $G$ is $\bigcap$-definable over $A$, then by Proposition 2.5, for any realization $a$ of $p$ there is some $Q$-internal $a_0 \in \mathrm{acl}(aA) - \mathrm{acl}(A)$. Replacing $a_0$ by the finite set of its $aA$-conjugates (which are all $Q$-internal as well), we may assume that there is an $A$-definable function $f$ with $f(a) = a_0$.

As the strong type $\mathrm{stp}(a/A)$ is a partial type over the finite set $Aa$, we may (after adding some imaginary parameters) assume that $p$ cosists of a single strong type. Put

$$H_0 = \{ g \in G : \exists x \models p \; [x \underset{A}{\downarrow} g \wedge xg \models p \wedge xg \underset{A}{\downarrow} g \wedge f(xg) = f(x)] \}.$$

Then $H_0$ is type-definable over $A$. Suppose $n, n' \in H_0$ with $n \underset{A}{\downarrow} n'$. The Independence Theorem yields the existence of $x \models p$ with $x \underset{A}{\downarrow} n, n'$, $xn \models p$ and $xn' \models p$, both independently from $n, n'$, such that $f(xn) = f(x) = f(xn')$. But since $xn \underset{A}{\downarrow} n, n'$ and $xn' \underset{A}{\downarrow} n, n'$, both $xn$ and $xn(n^{-1}n')$ realize $p$ independently over $A$ from $n^{-1}n'$, and $f(xn(n^{-1}n')) = f(xn') = f(xn)$. Hence $n^{-1}n' \in H_0$. By Lemma 1.10, $H_0 \cdot H_0 =: H$ is a subgroup of $G$ type-definable over $A$. So $H$ is $\bigcap$-definable by Lemma 2.3, and $H_0$ is generic in $H$.

Suppose $H$ has bounded index in $G$. Then $H_0$ is generic in $G$, so there is a generic $g \in H_0$ and independent $x \models p$ with $f(xg) = f(x)$. As $xg \underset{A}{\downarrow} x$, we get $f(xg) \underset{A}{\downarrow} xg$, contradicting $f(a) \notin \mathrm{acl}(A)$. Therefore $H$ has unbounded index in $G$, and by Lemma 2.3 it is contained in a definable supergroup $K$ of infinite index in $G$.

We claim that $G/K$ is $Q$-internal. In fact, since any two generic types are translates of one another (possibly by some independent new parameter), and any element of $G$ is the product of two generic elements, it is sufficient to check that $p/K$ is $Q$-internal. So consider $a \models p$. As it is generic for $G$ to stabilize $p$, there is an infinite generic Morley sequence $I := (x_i : i < \alpha)$ independent from $a$, such that $x_i a \models p$ for all $i < \alpha$. We claim that $aK$ is definable over $X := \{x_i, f(x_i a) : i < \alpha\}$ together with $A$. So suppose $\mathrm{tp}(a'/XA) = \mathrm{tp}(a/XA)$. Then $f(x_i a') = f(x_i a)$ for all $i < \alpha$, and $a' \underset{A}{\downarrow} I$. Now $a^{-1}a'$ must be independent over $A$ from $x_i a$ and from $x_i a'$ for some $i < \alpha$ (if $\alpha$ is big enough, since $(x_i a : i < \alpha)$ and $(x_i a' : i < \alpha)$ again form independent sequences), and we get $f(x_i a(a^{-1}a')) = f(x_i a') = f(x_i a)$. Therefore $a^{-1}a' \in H_0$, and $aK = a'K$. As $\mathrm{tp}(f(x_i a)/A)$ is $Q$-internal for all $i < \alpha$, the claim follows.

Finally, taking $\mathfrak{F}$ to be the family of all $G$-conjugates of $K$, Proposition 1.7 yields an $A$-definable normal subgroup $N$ of $G$ commensurable with a finite intersection of $G$-conjugates of $K$. Hence we may assume that $N$ contains a finite intersection of $G$-conjugates; since any $G$-conjugate of $G/K$ is again $Q$-internal, so is $G/N$. $\square$

It now follows that if $G$ is a $\bigcap$-definable group over a finite set $A$ in a small theory, and $Q$ is an $A$-invariant family of types, then there is a $Q$-connected $\bigcap$-definable component $G^Q$, which is unique up to commensurability. (This means that every generic type of $G^Q$ is foreign to $Q$, and for every definable supergroup $H$ of $G^Q$ the quotient $G/H$ is $Q$-analysable.)

**Question 2.7.** *If $Q$ is a family of formulas such that $G$ is $Q$-analysable, is there necessarily an analysis in finitely many steps?*

**Question 2.8.** *Do groups definable in a small simple theory satisfy property $\Re$ from [W2]? Does every such group have an infinite abelian, or finite-by-abelian subgroup?*

## 3. Supersimple Theories

In this section, $P$ will denote an $\emptyset$-invariant family of types closed under non-forking extensions, and $\mathrm{On}^+$ is the class of ordinals together with $\infty$ (where $\alpha < \infty$ for every ordinal $\alpha$). We shall quickly review the definition and basic properties of the $U_P$-rank defined in [W1] or [W2].

*Definition 3.1.* Let $q$ be a type. The $U_P$-*rank of $q$ relative to $P$*, denoted $U_P(q)$, is the smallest function from the collection of all types (over parameters in the monster model) to $\mathrm{On}^+$ satisfying for every ordinal $\alpha$:

$U_P(q) \geq \alpha + 1$ if there is an extension $q' \supseteq q$ over some set $A$, a type $p \in P$ over $A$, and realizations $a \models q'$ and $b \models p$ with $a \not\downarrow_A b$ and $U_P(a/Ab) \geq \alpha$.

Clearly, if $q$ is a nonforking extension of $p$, then $U_P(q) = U_P(p)$, and $U_P(q) = 0$ iff $q$ is hereditarily orthogonal to $P$. However, $q$ may well be foreign to $P$ and still have nonzero $U_P$-rank. In case $T$ is supersimple and $P$ is the family of all types, or of all types of $SU$-rank one, then $U_P$-rank is equal to $SU$-rank. By the usual proofs we obtain:

**Proposition 3.2.**
  (1) $U_P(a/bA) + U_P(b/A) \leq U_P(ab/A) \leq U_P(a/bA) \oplus U_P(b/A)$.
  (2) *If $a$ and $b$ are independent over $A$, then $U_P(ab/A) = U_P(a/A) \oplus U_P(b/A)$.*
  (3) *If $U_P(a/Ab) < \infty$ and $U_P(a/A) \geq U_P(a/Ab) \oplus \alpha$, then $U_P(b/A) \geq U_P(b/Aa) + \alpha$.*
  (4) *Suppose $U_P(a/Ab) < \infty$ and $U_P(a/A) \geq U_P(a/Ab) + \omega^\alpha \cdot n$. Then $U_P(b/A) \geq U_P(b/Aa) + \omega^\alpha \cdot n$.*

We shall in fact need a slight variant of 3.2(3).

**Lemma 3.3.** *If $U_P(a/AB) < \infty$ and $U_P(a/A) \geq U_P(a/AB) + 1$, then there is some $\bar{b} \in B$ with $U_P(\bar{b}/A) \geq U_P(\bar{b}/Aa) + 1$.*

*Proof.* By induction on $U_P(a/A)$, we may assume that the result holds for smaller $U_P$-rank. By definition of $U_P$-rank and the assumption of the lemma, there is some $C \supseteq A$ and $c$ realizing some type in $P$ over $C$ such that $a \not\downarrow_C c$ and $U_P(a/Cc) \geq U_P(a/AB)$. We may choose $Cc \downarrow_{Aa} B$, so $U_P(\bar{b}/Cca) = U_P(\bar{b}/Aa)$ for all $\bar{b} \in B$. Now if $B \not\downarrow_C c$, then there is some $\bar{b} \in B$ with $\bar{b} \not\downarrow_C c$, and

$$U_P(\bar{b}/A) \geq U_P(\bar{b}/C) \geq U_P(\bar{b}/Cc) + 1 \geq U_P(\bar{b}/Cca) + 1 = U_P(\bar{b}/Aa) + 1.$$

Otherwise $B \downarrow_C c$. Then since $a \not\downarrow_{BC} c$, we get

$$U_P(a/Cc) \geq U_P(a/AB) \geq U_P(a/CB) \geq U_P(a/CcB) + 1,$$

and by inductive hypothesis there is $\bar{b} \in B$ with

$$U_P(\bar{b}/A) \geq U_P(\bar{b}/Cc) \geq U_P(\bar{b}/Cca) + 1 = U_P(\bar{b}/Aa) + 1. \quad \square$$

This basically yields the local character of $P$-forking (which is obvious from the definition if $P$ is the family of all types):

**Corollary 3.4.** *If $A$ is the ascending union of sets $A_i$ for $i < \kappa$ and $U_P(a/A_i) = \alpha < \infty$ for all $i < \kappa$, then $U_P(a/A) = \alpha$.*

*Proof.* Suppose $U_P(a/A) < \alpha$. Then by Lemma 3.3 there is some $\bar{b} \in A$ with $U_P(\bar{b}/A_0) \geq U_P(\bar{b}/A_0 a) + 1$. By Proposition 3.2(3) we get $U_P(a/A_0) \geq U_P(a/A_0\bar{b}) + 1 \geq U_P(a/A_i) + 1$ for any sufficiently big $i$ such that $\bar{b} \in A_i$. But this contradicts our assumption. $\square$

We shall call a type $P$-*minimal* if every forking extension has smaller $U_P$-rank.

**Corollary 3.5.** *Every type of ordinal $U_P$-rank has a $P$-minimal extension of the same $U_P$-rank.*

*Proof.* As any chain of forking extensions has bounded length, this is immediate from Corollary 3.4. $\square$

**Corollary 3.5.** *Every $P$-minimal type $p$ of ordinal $U_P$-rank is non-orthogonal to a regular type.*

*Proof.* Let $q$ be a type of minimal $U_P$-rank non-orthogonal to $p$. Let $p$ be realized by $a$ and $q$ by $b$ over some set $A$; we may assume that $b$ is dominated by $a$ over $A$. We claim that $q$ is $P$-minimal. So suppose there is $c$ with $b \not\downarrow_A c$. Replacing $c$ by a finite bit of a Morley sequence in $\mathrm{tp}(b/Ac)$, we may assume that $\mathrm{tp}(c/A)$ has ordinal $U_P$-rank; we may clearly further assume $c \downarrow_{Ab} a$. BY domination $c \not\downarrow_A a$, whence $U_P(a/A) \geq U_P(a/Ac) + 1$ by $P$-minimality. By Proposition 3.2(3) we get

$$U_P(c/A) \geq U_P(c/Aa) + 1 \geq U_P(c/Aab) + 1 = U_P(c/Ab) + 1;$$

another application yields $U_P(b/A) \geq U_P(b/Ac) + 1$. This proves the claim.

Now $p$ is orthogonal to all types of $U_P$-rank less than $U_P(q)$, in particular to all forking extensions of $q$. As $a$ dominates $b$ over $A$, this shows that $q$ is orthogonal to all its forking extensions. $\square$

We shall also need a version of Shelah rank:

**Definition 3.6.** The *(Shelah) $D_P$-rank* is the least function from the class of all formulas to $\mathrm{On}^+$ satisfying

$D_P(\varphi) \geq \alpha + 1$ if there is some set $A$ of parameters containing those of $\varphi$, a realization $c$ of some type in $P$ over $A$, and a formula $\varphi'(x, c, A)$ forking over $A$, such that $\varphi' \subset \varphi$ and $D_P(\varphi') \geq \alpha$.

For a type $p$ we put $D_P(p) = \min\{D_P(\varphi) : \varphi \in p\}$. If $P$ is the family of all types, the subscript $P$ is omitted.

With this definition, $D_P(\varphi) = \max\{D_P(p) : p \in [\varphi]\}$. Clearly, if $\operatorname{tp}(c/a)$ is in $P$ and $a \not\!\!\downarrow_A c$, then $D_P(a/A) > D_P(a/Ac)$. In particular, if $q$ is a forking extension of $p$, then $D(q) < D(p)$. However, there is an example of a non-forking extension $q$ of a type $p$ with $D(q) < D(p)$. In general, $U_P(p) \leq D_P(p)$ for any type $p$.

## 4. Supersimple Groups

It is clear that to a type-definable group or quotient space we can ascribe a Lascar rank, namely the rank of any of its generic types. If a generic type of $G$ is $P$-minimal (and then all generic types are), then any type $p$ in $G$ is generic iff $U_P(p) = U_P(G)$, so a subgroup $H$ of $G$ has bounded index iff $U_P(H) = U_P(G)$; in any case the Lascar inequalities 3.2(1) immediately yield the corresponding inequalities for groups: if $H$ is a $\bigcap$-definable subgroup of $G$, then

$$U_P(H) + U_P(G/H) \leq U_P(G) \leq U_P(H) \oplus U_P(G/H).$$

Note however that although a type-definable group $G$ has Shelah rank, and in fact there is a generic type $p$ for $G$ with $D_P(G) = D_P(p)$, it is an open question whether there may be other generic types which have smaller Shelah rank.

**Lemma 4.1.** *Suppose $p$ is a $P$-minimal type in a simple group, with $U_P(p) = \omega^{\alpha_1} n_1 + \cdots + \omega^{\alpha_k} n_k =: \beta$, and for any two independent realizations $a$ and $b$ of $p$ we have $U_P(ab^{-1}) < \beta + \omega^{\alpha_k}$. Then $p$ is a translate of a generic type of the stabilizer $\operatorname{stab}(p)$. In particular the latter is $P$-minimal as well, and $U_P(\operatorname{stab}(p)) = \beta$.*

*Proof.* Recall that $\operatorname{stab}(p)$ contains a generic subset $S(p)$ such that $s \in S(p)$ iff there are realizations $a$ and $b$ of $p$, both independent from $s$, with $a = sb$, and $\operatorname{stab}(p) = S(p)^2$. It follows that $U_P(\operatorname{stab}(p)) \leq U_P(p)$.

Choose two independent realizations $a$ and $b$ of $p$. Since $U_P(ab^{-1}) \geq U_P(ab^{-1}/b) = U_P(a/b) = U_P(a) = \beta$, there is $A \downarrow_{ab^{-1}} a, b$ with $U_P(ab^{-1}/A) = \beta$, and such that $\operatorname{tp}(ab^{-1}/A)$ is $P$-minimal. Now $U_P(\bar{a}/a,b) = U_P(\bar{a}/ab^{-1})$ for all $\bar{a} \in A$; since $U_P(ab^{-1}) < U_P(ab^{-1}/\bar{a}) + \omega^{\alpha_k}$ by assumption, Proposition 3.2(3) yields first $U_P(\bar{a}) < U_P(\bar{a}/ab^{-1}) + \omega^{\alpha_k} = U_P(\bar{a}/a,b) + \omega^{\alpha_k}$, and then $U_P(a,b) < U_P(a,b/\bar{a}) + \omega^{\alpha_k}$. Since $a \downarrow b$ we have $U_P(a,b) = \beta \oplus \beta$; by Proposition 3.2(3) again $U_P(a,b) = U_P(a,b/\bar{a})$, and by $P$-minimality $a,b \downarrow \bar{a}$, whence $a,b \downarrow A$, and finally $a \downarrow_A b$. Now

$$U_P(ab^{-1}/Ab) = U_P(a/Ab) = \beta = U_P(ab^{-1}/A);$$

$P$-minimality of $\operatorname{tp}(ab^{-1}/A)$ yields $ab^{-1} \downarrow_A b$, and similarly $ab^{-1} \downarrow_A a$. By the definition of $S(p)$ we get $ab^{-1} \in S(p) \subseteq \operatorname{stab}(p)$. Hence $U_P(\operatorname{stab}(p)) = \beta$, and $\operatorname{tp}(ab^{-1}/A)$ is generic for $\operatorname{stab}(p)$ and a translate of $\operatorname{tp}(a/A) = p$. □

**Corollary 4.2.** *Suppose $U_P(G) = \omega^{\alpha_1} n_1 + \cdots + \omega^{\alpha_k} n_k$, and $\beta_i = \omega^{\alpha_1} n_1 + \cdots + \omega^{\alpha_i} n_i$ for some $1 \leq i \leq k$. Then $G$ has a type-definable normal subgroup $G_i$ of rank $\beta_i$; it is unique up to commensurability.*

*Proof.* If $p$ is a $P$-minimal type of $U_P$-rank $\beta_i$, it satisfies the requirements of Lemma 4.1. So $H := \operatorname{stab}(p)$ is a subgroup of $U_P$-rank $\beta_i$; if $H^g$ intersects $H$ in a subgroup of unbounded index, then $U_P(H \cap H^g) < \beta_i$ by $P$-minimality, whence $U_P(HH^g) \geq \beta_i + \omega^{\alpha_i}$ by the Lascar inequalities. This is greater than $U_P(G)$, contradiction. Therefore $H$ and any of its conjugates are commensurable, and the same holds for intersections of $H$ with its conjugates.

Therefore, if we fix a realization $g$ of some generic type $q$ of $G$ and let $\mathfrak{F}$ be the family of images of $H^g$ under strong automorphisms, then by Lemma 1.11 there is a type-definable commensurable subgroup $N$, which is invariant under automorphisms stabilizing $\mathfrak{F}$. But these include conjugation by $g^{-1}g'$ for independent realizations $g, g'$ of $q$, so $N$ is normalized by a subgroup of bounded index in $G$. Note that $N$ has only boundedly many conjugates over $\emptyset$, so after replacing it by the intersection of these conjugates, we may assume that it is type-definable over $\emptyset$. Let $S$ be a system of representatives for the cosets of $N_G(N)$ in $G$. Then the intersection $G_i := \bigcap_{s \in S} N^s$ is the required normal subgroup of $U_P$-rank $\beta_i$. □

Note that we cannot simply take $\bigcap_{g \in G} H^g$ as our normal subgroup commensurable with $H$: although it is normalized by $G$, since it is type-defined by $|G|$ formulas, we cannot be sure that it is normal in an elementary extension of $G$.

We shall now prove a definability result for groups type-definable in a supersimple theory.

**Proposition 4.3.** *Let $G$ be a type-definable group with $U(G) = \omega^\alpha n$ and $D(G) < \infty$. Then there is a definable supergroup $G_0$ of $G$, and definable subgroups $G_i$ of $G_0$ with $G = \bigcap_i G_i$.*

*Proof.* Clearly we may assume that all parameters are absorbed into the language, and that $G$ is $|T|^+$-saturated. First we show that there is a definable superset of $G$ on which multiplication is defined and associative, and inverses exist, and such that for any type containing $x \in X$ we have $U(p) \leq U(G)$. By compactness there is a definable superset $X_0$ of $G$ such that multiplication is defined and asociative on $X_0$ (but may go outside) and containing, for every $x \in X_0$, a (unique) inverse $x^{-1}$; we may assume it has minimal $D$-rank possible. Again by compactness, there is some definable $X_1 \subseteq X_0$ such that $X_1^2 \subseteq X_0$ and closed under inverse.

Suppose $X_1$ contains a type $p$ with $U(p) > U(G)$. Let $(x_i : i \in I)$ be a Morley sequence in $p$ and consider the sets $x_i G$, for $i \in I$. Since $X_1^2 \subseteq X_0$, they are contained in $X_0$. On the other hand, since $U(x_i/x_j) > U(G)$ for $i \neq j$, we get $x_i \notin x_j G$; by compactness and indiscernibility there is a definable superset $X'$ of $G$ contained in $X_1$ such that $x_i X' \cap x_j X' = \emptyset$ for all $i \neq j$. It follows that $D(X') < D(X_0)$, contradicting minimality of $D(X_0)$. Note that this does not use the fact that the $U$-rank of $G$ is a monomial. Let $X_2$ be a definable superset of $G$ closed under inverse, with $X_2^5 \subseteq X_1$.

Suppose $p$ is a type in $X_2$ with $U(p) = U(G)$. Then for any two independent realizations $a$ and $b$ of $p$ we get

$$U(G) \geq U(X_1) \geq U(ab^{-1}) \geq U(ab^{-1}/b) = U(a/b) = U(a) = U(G),$$

so $ab^{-1}$ is independent from $b$ and similarly from $a$. Therefore $H := \text{stab}(p)$ is a type-definable group contained in $X_2^4$ with $U(H) = U(p)$; it must be commensurable with $G$ since otherwise $U(GH) > U(G)$, but $GH \subseteq X_2 X_2^4 \subseteq X_1$, contradicting $U(X_1) = U(G)$. Hence $D(H, \varphi, k) = D(G, \varphi, k)$ for all stratified formulas $\varphi$ and all $k < \omega$; since $p$ is a translate of a generic type for $H$, we obtain $D(p,\varphi,k) = D(G,\varphi,k)$ for all such $\varphi, k$.

Conversely, for any enumeration $(\varphi_i, k_i)$ of stratified formulas and natural numbers, suppose $p$ is a type in $X_2$ with $D(p,\varphi_i,k_i) = n_i$ is maximal possible subject to $D(p,\varphi_j,k_j) \geq n_j$ for all $j < i$. Let $a$ realize $p$ and $b$ realize a generic type for $G$ independently from $a$. Then

$$U(G) = U(X_1) \geq U(ab) \geq U(ab/a) = U(b/a) = U(b) = U(G),$$

so we must have equality and $ab$ is independent from $a$. Therefore

$$\begin{aligned} D(b,\varphi_i,k_i) = D(b/a,\varphi_i,k_i) &= D(ab/a,\varphi_i,k_i) \\ &= D(ab,\varphi_i,k_i) \geq D(ab/b,\varphi_i,k_i) = D(a/b,\varphi_i,k_i) = D(a,\varphi_i,k_i) = n_i \end{aligned}$$

for all $i$; since $n_i$ was chosen maximal, this implies that we have equality all the way through and $ab$ is independent from $b$. Therefore $U(a) = U(a/b) = U(ab/b) = U(ab) = U(G)$.

It follows that the set $P$ of types over $G$ in $X_2$ of maximal $U$-rank is closed as the set of those types $p \in S_1(G)$ containing the formula $x \in X_2$ and such that $D(p,\varphi_i,k_i) = n_i$ for all $i$. Choose some definable superset $X_3$ of $G$ with $X_3^2 \subseteq X_2$.

For a definable superset $Y$ of $G$ contained in $X_3$, consider the relation $R(x,y)$ iff there are independent realizations $x', y'$ of the same Lascar strong type of $x$ and $y$ over $G$, such that $x'Y \cap y'Y$ does not fork over $G$. Suppose $a$ and $b$ independently realize types in $P$ such that $aY \cap bY$ does fork over $G$, and let $(a_i, b_i : i \leq n)$ be a Morley sequence in $\text{tp}(a,b/G)$. If $c \in \bigcap_{i \leq n} a_i Y \cap b_i Y$, then $c$ must fork with every $a_i, b_i$ over $G \cup (a_j, b_j : j < i)$ and this must affect the coefficient of $\omega^\alpha$ in the $U$-rank of $\text{tp}(c/G, a_j, b_j : j < i)$; as it can go down at most $n$ many times, there can be no such $c$ and the intersection is empty. Therefore $\neg R(x,y)$ is a closed condition on $P^2$; since $R$ is obviously closed as well, by compactness there is a formula $\varphi(x,y)$ with parameters in $G$ which agrees with $R(x,y)$ on $P \times P$. As $R$ is a stable relation by [KP], there also is a formula $P_0$ in the partial type $P$ and some $n < \omega$ such that $\rho(x,y) := x \in P_0 \land y \in P_0 \land \varphi(x,y)$ can order sequences of length at most $n$. In other words, $\rho$ is a stable formula. It follows that the set

$$\psi(y) := \text{for all generic } g \text{ we have } \rho(g,y)$$

is definable, say over some parameters $a \in G$. Note that for $y \models P$ the formula $\psi(y)$ holds iff for all generic $g$ the set $gY \cap yY$ does not fork over $G$.

Suppose $y \in X_3$ and $g, g' \in G$ are two independent elements generic over $a$ and independent from $y$ over $a$. Then $g \downarrow_{a,y} g'g^{-1}$, so $gy \downarrow_{a,y} g'g^{-1}$; as $y \downarrow_a g'g^{-1}$ we get $gy \downarrow_a g'g^{-1}$, and $U(gy/a, g'g^{-1}) = U(G)$. So $\psi(gy)$ holds iff $\psi(h)$ holds for any realization $h$ of the non-forking extension of $\text{tp}(gy/a, g'g^{-1})$ to $G$ iff for all generic $g''$ the set $g''Y \cap hY$ does not fork over $G$. But if $g''$ runs through realizations for all generic types of $G$ over $G$, so does $g'g^{-1}g''$. Therefore $g''Y \cap hY$ does not fork over $G$ for that $h$ and all generic $g''$ iff $g'g^{-1}g''Y \cap g'g^{-1}hY$ does not fork over $G$ for all generic $g''$ iff $g''Y \cap g'g^{-1}hY$ does not fork over $G$ for all generic $g''$ iff $\models \psi(g'g^{-1}h)$ iff $\models \psi(g'g^{-1}gy)$ iff $\models \psi(g'y)$. In other words, if $\psi(gy)$ holds for some generic $g \in G$ independent from $y$ over $a$, it holds for all generic $g' \in G$ such that $g, g'$ is independent from $y$ over $a$.

Let $g_0, \ldots, g_{2n}$ be $2n+1$ elements of $G$ independent and generic over $a$. For any $y \in X_3$ there are at least $n+1$ of the $g_i$ which are (collectively) independent from $y$ over $a$, as the coefficient of $\omega^\alpha$ in the $U$-rank must go down every time $y$ forks with some $g_i$ over some others. Therefore

$$\psi_Y(y) = \bigvee_{i_0,\ldots,i_n \leq 2n} \bigwedge_{k=0}^{n} \psi(g_{i_k} y)$$

holds of $y \in X_3$ iff $\models \psi(gy)$ for some $g \in G$ generic over $a$ and independent from $y$ over $a$ iff $\models \psi(gy)$ for such $g \in G$ iff $g'Y \cap gyY$ does not fork over $G$ for all generic $g'$ (over $G \cup \{y\}$) and all $g \in G$ generic and independent from $y$ over $a$.

Suppose $y \in X_3$ satisfies $\psi_Y$ for all $\emptyset$-definable supersets $Y$ of $G$. Let $g' \in G$ be generic over $y$ (and the parameters used for all $\psi_Y$) and $g$ be generic over $G$. Then for any such $Y$ the intersection $gY \cap g'yY$ is non-empty, and there are elements $h, h'$ in (some elementary extension of) $G$ with $gh = g'yh'$, whence $y \models G$. By saturation of $G$, this shows that $\bigcap \psi_Y$ type-defines a subset of $G$. On the other hand, if $y \models G$ and $g, g'$ are any elements in $G$, then $gY \cap g'yY$ contains $G$ and thereby elements which do not fork over $G$, whence any element of $G$ satisfies $\psi_Y$. Thus $\bigcap \psi_Y$ type-defines $G$.

Now suppose $g \in G$, $y \in X_3$, $gy \in X_3$ and $y \models \psi_Y$ for some definable superset $Y$ of $G$. If $h \in G$ is generic over $g, y$ and the parameters used for $\psi_Y$, so is $hg$, and $h$ is generic over $gy$. Therefore $\models \psi_Y(y)$ iff $\models \psi(hy)$ iff $\models \psi((hg)y)$ iff $\models \psi(h(gy))$ iff $\models \psi_Y(gy)$. So for $Y$ small enough such that $\psi_Y^2 \subseteq X_3$, the set

$$S_Y := \{x \models \psi_Y : \forall y \, [\psi_Y(y) \to \psi_Y(xy)]\}$$

defines a superset of $G$ closed under multiplication, and the invertible elements of $S_Y$ form a definable supergroup $G_Y$ of $G$. Clearly, the intersection of all $G_Y$ is equal to $G$. □

**Theorem 4.4.** *A type-definable group $G$ in a supersimple theory is $\bigcap$-definable.*

*Proof.* Clearly, we may assume that the ambient model is sufficiently saturated. We use induction on the number of monomials in $U(G)$. So suppose $U(G) = \omega^\alpha n + \beta$ with $\beta < \omega^\alpha$. By Corollary 4.2 there is a type-definable normal subgroup $N$ of $G$ with $U(N) = \omega^\alpha n$, and $N$ is $\bigcap$-definable by Proposition 4.3. So if $X$ is a definable superset of $G$ on which multiplication is defined, associative, and inverses exist, then there is a sequence of definable groups $H_i$ contained in $X$ with $\bigcap_i H_i = N$; by supersimplicity we may assume $U(H_i) = U(N)$ for all $i$. By compactness there is some $i$ such that $(GH_i)^4 \subseteq X$. As $U(H_i^g) = U(N) \leq U(H_i^g \cap H_i^{g'})$ for all $g, g' \in G$, this means that the intersection of any two $G$-conjugates of $H_i$ has finite index in either conjugate; by saturation there is a bound on this index, and by Fact 1.5 (since the product of any two $G$-conjugates of $H_i$ is contained in $X$) there is a $G$-invariant definable group $H$ commensurable with $H_i$, such that $GH \subseteq X$. Then $N_X(H)/H$ is interpretable, and so is $GH/H$; as $U(GH/H) = \beta$, the inductive assumption yields $GH/H$ to be $\bigcap$-definable. In particular there is a definable group $G_0 \subseteq N_X(H)$ with $G/H \leq G_0/H$, whence $G \leq G_0 H \subseteq X$. Therefore every definable superset of $G$ contains a definable supergroup of $G$; it follows that $G$ is $\bigcap$-definable. □

**Corollary 4.5.** *A type-definable division ring $F$ in a supersimple theory is definable.*

*Proof.* By Theorem 4.4 the additive and the multiplicative group of $F$ are both $\bigcap$-definable. Let $M$ be a multiplicative supergroup of $F$ and $A, A'$ additive supergroups of $F$ with $A \supseteq M \cup \{0\} \supseteq A' \supseteq F$ and $U(A) = U(F)$; by compactness we may choose these small enough such that the distributive law holds, and $kA' + k'A' \subseteq M$ for all $k, k' \in F^\times$. Then for every $k \in F^\times$ the group $kA'$ is an isomorphic image of $A'$ inside $A$; since $U(A) = U(A')$ this means that $A' \cap kA'$ has finite index in $A'$ and in $kA'$. By compactness this index is bounded; by Fact 1.5 there is a definable $F^\times$-invariant subgroup $A_0$ commensurable with $A'$ contained in $kA' + k'A'$ for some $k, k' \in F^\times$, and hence $A_0 \subseteq M$. Let $F_0 := \{m \in A_0 : mA_0 \leq A_0\}$, a definable subring of $A_0$ containing $F$. Clearly $F_0$ has no zero divisors. But then for any $0 \neq a \in F_0$ the sequence $(a^i F_0 : i < \omega)$ cannot descend infinitely, as otherwise the formula $xA_0 < yA_0$ would define a partial order with infinite chains, contradicting simplicity. Hence $a^i F_0 = a^{i+j} F_0$ for some $j > 0$ and $a^i(1 - a^j k) = 0$ for some $k \in F_0$. It follows that $a$ is invertible, and $F_0$ is a division ring.

Starting with a smaller $A$, we see that $F$ is the intersection of definable division rings $F_i$; as the additive index $|F_i^+ : F_j^+|$ is infinite for any $i < j$ and therefore $U(F_i) > U(F_j)$, supersimplicity yields $F_i = F_j$, whence $F = F_0$. $\square$

The group decomposition can beused to prove that for a generic type of a supersimple goup Lascar strong type is equal to strong type.

**Lemma 4.6.** *Let $H$ be an $A$-definable coset of a group in a simple theory with $U(H) = \omega^\alpha n$, and suppose $p$ is a generic type for $H$ over $A$. Let $q$ be a type over $A$ such that for any realization $a$ of $q$ there is a realization $a'$ of $p$ such that $a$ and $a'$ are interalgebraic over $A$. Then the equivalence relation $r(x, y)$, which holds for two realizations of $q$ iff they have the same Lascar strong type, is given by an interscetion of finite $A$-definable equivalence relations.*

*Proof.* Recall that $r$ is type-definable by [K3]. Let $\varphi(x, y)$ be any formula in $r$, and define the relation $R_\varphi(a, b)$ to hold iff for some independent $a'$ and $b'$ of the same Lascar strong types as $a$ and $b$ over $A$, the formula $\varphi(a', z) \wedge \varphi(z, b')$ does not fork over $A$. This is clearly invariant under Lascar strong type; note that if $b''$ is independent from $a$ over $A$ and has the same Lascar strong type as $b'$ and $\varphi(a', z) \wedge \varphi(z, b')$ does not fork over $A$, then there is $c \downarrow_A a', b'$ realizing $\varphi(a', z) \wedge \varphi(z, b')$, and by the Independence Theorem we can find $c' \downarrow_A a'b''$ with $\operatorname{tp}(c'a'/A) = \operatorname{tp}(ca'/A)$ and $\operatorname{tp}(c'b''/A) = \operatorname{tp}(cb'/A)$. Hence whether or not $\varphi(a', z) \wedge \varphi(z, b')$ forks over $A$ does not depend on the particular choice of $b'$ (or $a'$).

Clearly $R_\varphi$ is a closed relation. Let $\vartheta(x, y)$ be the formula witnessing interalgebraicity of realizations of $p(x)$ and $q(y)$, and let $X := \{y : \exists x \in H \; \vartheta(x, y)\}$. Then any element of $X$ has Lascar rank at most $\omega^\alpha n$. Since forking with some element in a Morley sequence in $\operatorname{stp}(a'b'/A)$ must affect the coefficient of $\omega^\alpha$, if $\varphi(a', z) \wedge \varphi(z, b')$ does fork over $A$, then for any Morley sequence $(a_i, b_i : i \leq n)$ in $\operatorname{stp}(a', b'/A)$ the formula $\bigwedge_{i \leq n} \varphi(a_i, z) \wedge \varphi(z, b_i)$ is inconsistent with $X$. Hence $\neg R_\varphi$ is closed as well, and $R_\varphi$ is a definable relation.

Choose any realization $a$ of $q$, and define $E_a(x, y)$ as $R_\varphi(a, x) \leftrightarrow R_\varphi(a, y)$. Then $E$ is a definable equivalence relation invariant under Lascar strong type. So $E_a$ and any of its $A$-conjugates is coarser than equality of Lascar strong type; compactness now yields that the intersection of all $A$-conjugates of $E_a$ is a finite intersection, and hence a finite $A$-definable equivalence relation $E_\varphi$. (Note that $E_\varphi$ will be a finite equivalence relation on some definable superset $\psi$ of $q$; we can then declare $\neg\psi$ to be a single equivalence class and have $E_\varphi$ defined on the whole of $X$.)

We finally claim that $r$ is the intersection of all the $E_\varphi$, for $\varphi \in r$. Indeed, if $r(a,b)$ holds, then since $R_\varphi$ is invariant under Lascar strong type, for any $c \models p$ we have $R(c,a)$ iff $R(c,b)$. Thus $E_\varphi(a,b)$ holds for all $\varphi \in r$. Conversely, suppose $E_\varphi(a,b)$ holds for all $\varphi \in r$. Since $r$ is transitive, for any $\varphi \in r$ there is a $\psi \in r$ such that $\exists zz'\, \psi(x,z) \wedge \psi(z,z') \wedge \psi(z',y$ implies $\varphi(x,y)$. Since $E_\psi(a,b)$ holds, we have $R_\psi(c,a)$ iff $R_\psi(c,b)$ for any $c \models q$, in particular for $c = a$. But clearly $R_\psi(a,a)$ holds, whence $\psi(a,z) \wedge \psi(z,b')$ does not fork over $A$ for any $b'$ of the same Lascar strong type as $b$ independent from $a$ over $A$. In particular this formula is consistent; as in addition $\models \psi(b',b)$, we get $\models \varphi(a,b)$. □

**Theorem 4.7.** *Let $p$ be the generic type of a type-definable group in a supersimple theory. Then Lascar strong type is equal to strong type for $p$.*

*Proof.* Suppose $p$ is over $A$. We may assume that $A$ is algebraically closed; it is sufficient to show that if $B$ and $B'$ are two sets independent over $A$ and $q$ and $q'$ are two non-forking extensions of $p$ to $C$ and $C'$ respectively, then there is $a \downarrow_A CC'$ realizing both $q$ and $q'$.

Note that the group $G$ for which $p$ is generic is definable over $A$. We shall use induction on $U(G)$ and prove the theorem more generally for types which are interalgebraic with $p$. If $U(G) = \omega^\alpha n + \beta$ with $\beta < \omega^\alpha$, then by Corollary 4.2 there is an $A$-definable normal subgroup $N$ of $G$ with $U(N) = \omega^\alpha n$. Let $a \models q$ and $a' \models q'$. Then $aN$ and $a'N$ are generic types for $G/N$ of the same type over $A$ and rank $\beta$, and independent from $C$ and $C'$ respectively. By induction, we may assume $aN = a'N$, and even $\mathrm{acl}(aN,A) = \mathrm{acl}(a'N,A) = A'$, whence $a$ and $a'$ have the same type over $A'$. (It is here that we have to consider sequence of elements interalgebraic with $aN$ and $a'N$, rather than simply the elements themselves.) But now $aN \downarrow_A CC'$, whence $C \downarrow_{A'} C'$, and clearly $a \downarrow_{A'} C$ and $a' \downarrow_{A'} C'$. Furthermore $a$ and $a'$ realize a (the same) generic type of the coset $aN$, so if $\bar{a}$ is interalgebraic with $a$ and $\bar{a}'$ is interalgebraic with $a'$ such that $\mathrm{tp}(\bar{a}/A) = \mathrm{tp}(\bar{a}'/A)$, then $\mathrm{tp}(bara/A') = \mathrm{tp}(bara'/A')$, and by Lemma 4.6 there is $\bar{c} \downarrow_{A'} CC'$ with $\mathrm{tp}(\bar{c}/CA') = \mathrm{tp}(\bar{a}/CA')$ and $\mathrm{tp}(\bar{c}/C'A') = \mathrm{tp}(\bar{a}'/C'A')$. Finally $CC' \downarrow_A A'$ implies $\bar{c} \downarrow_A CC'$, so the theorem is proved. □

We shall now formulate the appropriate version of the Zil'ber Indecomposability Theorem. (Note that we revert to the case of an arbitrary family $P$; if the theory is supersimple, we shall get definable groups by Theorem 4.4.)

**Theorem 4.8.** *Suppose $G$ is a type-definable group in a simple theory of rank $U_P(G) = \omega^\alpha n + \beta$ with $\beta < \omega^\alpha$, and $\{X_i : i \in I\}$ is a family of type-definable sets. Then there is a type-definable subgroup $H$ of $G$ such that for any $i \in I$ we have $U_P(X_iH) < U_P(H) + \omega^\alpha$ (meaning that $U_P(p) < U_P(H) + \omega^\alpha$ for any type $p$ containing the formula $x \in X_iH$). Furthermore, there are $i_1,\ldots,i_m \in I$ such that $H \subseteq X_{i_1}^{\pm 1} \cdots X_{i_m}^{\pm 1}$.*

Note that if $H$ is $\bigcap$-definable, this translates as $U_P(X_i/H) < \omega^\alpha$ (in fact, this also makes sense for type-definable $H$ by using hyperimaginaries). So if $X_i$ is $\alpha_P$-indecomposable and contains 1, then $X_i$ is contained in $H$; if this holds for all $X_i$, then $H$ is the group generated by $\{X_i : i \in I\}$.

*Proof.* Let $k$ be maximal such that there are $i_1, \ldots, i_m$, such that $X_{i_1}^{\pm 1} \cdots X_{i_m}^{\pm 1} =: Y$ contains a type of $U_P$-rank at least $\omega^\alpha k$. Choose a $P$-minimal type $p$ containing the formula $x \in Y$ with $U_P(p) = \omega^\alpha k$, and put $H := \mathrm{stab}(p)$. Note that $S(p) \subseteq Y^2$, and $\mathrm{stab}(p) \subseteq Y^4$. By maximality of $k$, the assumptions of Lemma 4.1 are satisfied, and $U_P(H) = \omega^\alpha k$. Clearly $U_P(X_i H) < U_P(H) + \omega^\alpha$ for any $i \in I$; in fact we even have $U_P(XH) < U_P(H) + \omega^\alpha$ for any finite product $X$ of elements of the family $\{X_i : i \in I\}$. $\square$

Note that if $H'$ is a conjugate of $H$ under an automorphism stabilizing $Y$ setwise, then by $P$-minimality either $U_P(H \cap H') < \omega^\alpha k$, or $H$ and $H'$ are commensurable. But the first case cannot happen, since then $U_P(HH') \geq \omega^\alpha(k+1)$, contradicting maximality of $k$. Therefore as above in the proof of Corollary 4.2 we can apply Lemma 1.11 to find a subgroup type-definable over the parameters used for $Y$ (and still contained in $Y$) which is commensurable with $H$; furthermore, if all the $X_i$ are $G$-invariant subsets, we may obtain a normal subgroup for $H$.

In particular, if $T$ is small or superstable and $G$ and all $X_i$ are type-defined over a finite set, then $H$ can be chosen type-defined over a finite set and $\bigcap$-definable by Lemma 2.3 or Theorem 4.4. If in addition all the $X_i$ are definable, then we actually get a big definable subgroup inside $Y^4$ by compactness.

Recall that a group is simple if it has no normal subgroups, and definably simple if it has no definable normal subgroups.

**Corollary 4.9.** *Suppose $G$ is a type-definable, definably simple group in a simple theory with $0 < U_P(G) < \infty$. Then $G$ is simple.*

*Proof.* Suppose $U_P(G) = \omega^\alpha n + \beta$ with $\beta < \omega^\alpha$. We claim first that $U_P(g^G) \geq \omega^\alpha$ for any $g \neq 1$ in $G$. So suppose otherwise. Then $U_P(C_G(g)) \geq \omega^\alpha n$ by the Lascar inequality for groups (since $U_P(G/C_G(g)) = U_P(g^G)$). By Proposition 1.7 there is a definable normal subgroup $N$ commensurable with a sufficiently big finite intersection of conjugates of $C_G(g)$. As $U_P(G/C_G(g)) < \omega^\alpha$ and the same holds for any conjugate, we see that $U_P(G/N) < \omega^\alpha$, and $N$ is non-trivial. On the other hand, $C_G(g)$ cannot have finite index in $G$, as otherwise the intersection of its conjugates would form a definable subgroup of finite index; since also $Z(G)$ is trivial, this index must be greater than 1 and we obtain a contradiction.

Next we claim that $G$ has no type-definable normal subgroup. For if $N$ is type-definable normal minimal (up to commensurability) and $n \in N$ is nontrivial, then $X := (n^G)^k$ contains a type-definable normal subgroup $N_1$ for big $k < \omega$ by Theorem 4.8; by minimality the index of $N_1$ in $N$ must be bounded. But this means that finitely many translates of $X$ must cover $N$, and the subunion of those translates by elements in $N$ must be equal to $N$. This contradicts definable simplicity. (Note that for supersimple theories, definable and type-definable simplicity are easily seen to be the same by Theorem 4.4 and Fact 1.5; the argument above is need only for the case of a more general family $P$.)

Finally, if $N$ were a non-trivial normal subgroup of $G$ and $1 \neq n \in N$, then $(n^G)^k$ contains a non-trivial type-definable normal subgroup of $G$ for big $k < \omega$, which must be proper since it is contained in $N$. This final contradiction finishes the proof. □

**Question 4.9.** *If $U_P(G) = 1$, does this imply that $G$ has a finite-by-abelian subgroup $A$ with $U_P(A) = 1$? If $D$ is a supersimple division ring, is it commutative? (The answer is yes if $D$ contains enough roots of unity.) Is a supersimple field bounded pseudo-algebraically closed?*

## 5. One-based Groups

**Definition 5.1.** A structure is *one-based* if for any two tuples $a$ and $b$ are independent over $\operatorname{acl}(a) \cap \operatorname{acl}(b)$.

Note that if a structure is one-based, then every type $\operatorname{tp}(a/A)$ has a canonical base $\operatorname{acl}(a) \cap \operatorname{acl}(A)$.

**Proposition 5.2.** *Let $G$ be a one-based group $\emptyset$-definable in a simple theory. Then every definable subgroup $H$ is commensurable with one definable over $\operatorname{acl}(\emptyset)$.*

*Proof.* Let $H^c$ be the locally connected component as given by Corollary 1.8, with canonical parameter $u$. Let $h$ realize a generic type for $H^c$ over $u$, and $g$ realize a generic type for $G$ over $u, h$. So $\operatorname{tp}(hg/g, u)$ is generic for $H^c g$ over $u, g$, and $\operatorname{tp}(hg/u, h)$ is generic for $G$.

Now let $v$ be the canonical parameter for the coset $H^c g$. If $(h_i : i < \omega)$ is a Morley sequence in $\operatorname{tp}(hg/g, u)$, then $h_i$ is generic in $H^c g$ for all $i < \omega$. But any distinct coset of a conjugate of $H^c$ intersects $H^c g$ in a coset of infinite index and cannot contain infinitely many generic elements, so $v$ is definable over $(h_i : i < \omega)$. As $v \in \operatorname{dcl}(g, u)$, we get $v \in \operatorname{Cb}(hg/g, u)$. On the other hand, one-basedness implies $\operatorname{Cb}(hg/g, u) \subset \operatorname{acl}(hg)$.

As $H^c = (H^c g)(H^c g)^{-1}$, the canonical parameter $u$ is definable over $v$, whence $u \in \operatorname{acl}(hg)$. But $hg \underset{\bigsmile}{} u$, so $u \in \operatorname{acl}(\emptyset)$, and $H^c$ is definable over $\operatorname{acl}(\emptyset)$. □

**Corollary 5.3.** *If $\mathfrak{F}$ is a family of uniformly definable subgroups of a one-based group definable in a simple theory, then there are only finitely many commensurability classes among members of $\mathfrak{F}$.*

*Proof.* Assuming the Corollary does not hold, compactness yields a subgroup which is not commensurable to any $\operatorname{acl}(\emptyset)$-definable subgroup, contradicting Proposition 5.2. □

**Theorem 5.4.** *Suppose $G$ is a one-based group definable in a simple theory. Then $G$ is finite-by-abelian-by-finite.*

*Proof.* Consider, for $g \in G$, the subgroup $H_g = \{(h, h^g) : h \in G\} < G^2$. Then $\mathfrak{F} := \{H_g : g \in G\}$ is a family of uniformly definable subgroups of $G^2$ and has only finitely many commensurability classes. But $H_g$ and $H_{g'}$ are commensurable iff $C_G(g'g^{-1})$ has finite index in $G$. If we define $Z^* := \{g \in G : |G : C_G(g)| < \omega\}$, then the commensurability classes of $\mathfrak{F}$ correspond to the cosets of $Z^*$ in $G$. Therefore $Z^*$ has finite index in $G$; replacing $G$ by $Z^*$ we may assume that every centralizer has finite index in $G$. But then $[g, G]$ is finite for every $g \in G$, and for independent $g, h$ in $G$ we get
$$[g, h] \in \mathrm{acl}(g) \cap \mathrm{acl}(h) = \mathrm{acl}(\emptyset).$$

But since $C_G(g)$ has finite index in $G$ for every $g$, for every $h$ there is $h'$ independent from $g$ with $[g, h] = [g, h']$. It follows that there are only finitely many commutators, and $G'$ is finite. $\square$

In fact, replacing $G$ by $C_G(G')$, we may even assume that $G'$ is central in $G$.

**Question 5.5.** *Is there a characterization of one-based simple groups in terms of their definable subsets?*

MATHEMATICAL INSTITUTE, 24–29 ST GILES', OXFORD OX1 3LB, ENGLAND
THE FIELDS INSTITUTE, 222 COLLEGE STREET, TORONTO, ONTARIO, CANADA K5T 3J1